\documentclass{amsart}

\usepackage{amsmath}
\usepackage{amscd}
\usepackage{amssymb}

\newcommand{\cal}{\mathcal}

\newcommand{\bC}{{\Bbb C}}
\newcommand{\bE}{{\Bbb E}}
\newcommand{\bL}{{\Bbb L}}

\newcommand{\bP}{{\Bbb P}}
\newcommand{\bQ}{{\Bbb Q}}

\newcommand{\bZ}{{\Bbb Z}}

\newcommand{\cC}{{\cal C}}

\newcommand{\cM}{{\cal M}}
\newcommand{\cO}{{\cal O}}

\newcommand{\Mbar}{\overline{\cM}}

\DeclareMathOperator{\Aut}{Aut}
\DeclareMathOperator{\End}{End}

\DeclareMathOperator{\id}{id}

\DeclareMathOperator{\tr}{tr}

\newtheorem{theorem}{Theorem}[section]
\newtheorem{theorem/definition}{Theorem/Definition}[section]
\newtheorem{Theorem}{Theorem}

\newtheorem{lemma}{Lemma}[section]

\theoremstyle{remark}

\theoremstyle{definition}

\begin{document}

\title{On a proof of a conjecture of Marino-Vafa on Hodge Integrals}
\author{Chiu-Chu Melissa Liu}
\address{Department of Mathematics\\
Harvard University\\ Cambridge, MA 02138, USA}
\email{ccliu@math.harvard.edu}
\author{Kefeng Liu}
\address{Department of Mathematics\\
UCLA\\ Los Angeles, CA 90095-1555}
\email{liu@math.ucla.edu}
\author{Jian Zhou}
\address{Department of Mathematical Sciences\\
Tsinghua University\\ Beijing, 100084, China}
\email{jzhou@math.tsinghua.edu.cn}
\begin{abstract}
We outline a proof of a remarkable formula for Hodge integrals conjectured by
Mari\~{n}o and Vafa \cite{Mar-Vaf} based on large $N$ duality.
\end{abstract}

\maketitle

\section{Introduction}

Let $\Mbar_{g,n}$ denote the Deligne-Mumford moduli stack of stable curves of genus
$g$ with $n$ marked points. Let $\pi:\Mbar_{g,n+1}\to \Mbar_{g,n}$
be the universal curve, and let $\omega_\pi$ be the relative dualizing sheaf.
The Hodge bundle
$$
\bE=\pi_*\omega_\pi
$$
is a rank $g$ vector bundle over $\Mbar_{g,n}$ whose
fiber of over $[C,x_1,\ldots,x_n]\in\Mbar_{g,n}$ is
$H^0(C,\omega_C)$.
Let $s_i:\Mbar_{g,n}\to\Mbar_{g,n+1}$ denote the section of $\pi$ which
corresponds to the $i$-th marked point, and let
$$
\bL_i=s_i^*\omega_\pi
$$
be the line bundle over $\Mbar_{g,n}$ whose fiber over
$[C,x_1,\ldots,x_n]\in\Mbar_{g,n}$ is the cotangent line $T^*_{x_i} C$
at the $i$-th marked point $x_i$.
A Hodge integral is an integral of the form
$$\int_{\Mbar_{g, n}} \psi_1^{j_1}
\cdots \psi_n^{j_n}\lambda_1^{k_1} \cdots \lambda_g^{k_g}$$
where
$\psi_i=c_1(\bL_i)$ is the first Chern class of $\bL_i$, and
$\lambda_j=c_j(\bE)$ is the $j$-th Chern class of the Hodge bundle.

Hodge integrals arise in the calculations of Gromov-Witten invariants
by localization techniques \cite{Kon2, Gra-Pan}.
The explicit evaluation of Hodge integrals is a difficult problem.
The Hodge integrals involving only $\psi$ classes
can be computed recursively by Witten's conjecture \cite{Wit}
proven by Kontsevich \cite{Kon1}. Algorithms of computing Hodge
integrals are described in \cite{Fab}.

In \cite{Mar-Vaf}, M. Mari\~{n}o and C. Vafa obtained a closed
formula for a generating function of certain open Gromov-Witten
invariants, some of which has been reduced to Hodge integrals
by localization techniques which are not fully clarified
mathematically.
This leads to a conjectural formula of Hodge integrals. To state this
formula, we introduce some notation, following \cite{Zho3}.
Let
$$
\Lambda^{\vee}_g(u)=u^g-\lambda_1 u +\cdots+(-1)^g\lambda_g
$$
be the Chern polynomial of $\bE^\vee$, the dual of the
Hodge bundle. For a partition $\mu$ given by
$$
\mu_1\geq\mu_2\geq\cdots\geq\mu_{l(\mu)} >0,
$$
let $|\mu|=\sum_{i=1}^{l(\mu)}\mu_i$, and define
\begin{eqnarray*}
\cC_{g, \mu}(\lambda; \tau)& = & - \frac{\sqrt{-1}^{|\mu|+l(\mu)}}{|\Aut(\mu)|}
 [\tau(\tau+1)]^{l(\mu)-1}
\prod_{i=1}^{l(\mu)}\frac{ \prod_{a=1}^{\mu_i-1} (\mu_i \tau+a)}{(\mu_i-1)!} \\
&& \cdot \int_{\Mbar_{g, l(\mu)}}
\frac{\Lambda^{\vee}_g(1)\Lambda^{\vee}_g(-\tau-1)\Lambda_g^{\vee}(\tau)}
{\prod_{i=1}^{l(\mu)}(1- \mu_i \psi_i)}, \\
\cC_{\mu}(\lambda; \tau) & = & \sum_{g \geq 0} \lambda^{2g-2+l(\mu)}\cC_{g,
\mu}(\lambda; \tau)
\end{eqnarray*}

Note that
$$
\int_{\Mbar_{0, l(\mu)}}
\frac{\Lambda^\vee_0(1)\Lambda^\vee_0(-\tau-1)\Lambda_0^\vee(\tau) }
{\prod_{i=1}^{l(\mu)} (1 -\mu_i \psi_i) }
=\int_{\Mbar_{0, l(\mu)}}
\frac{1}{\prod_{i=1}^{l(\mu)}(1 - \mu_i\psi_i)}
= |\mu|^{l(\mu)-3}
$$
for $l(\mu)\geq 3$, and we use this expression to extend the definition
to the case $l(\mu)<3$.

Introduce formal variables $p=(p_1,p_2,\ldots,p_n,\ldots)$, and define
$$
p_\mu=p_{\mu_1}\cdots p_{\mu_{l(\mu)} }
$$
for a partition $\mu=(\mu_1\geq \cdots \geq \mu_{l(\mu)}>0 )$. Define
generating functions
\begin{eqnarray*}
\cC(\lambda; \tau; p) & = & \sum_{|\mu| \geq 1} \cC_{\mu}(\lambda;\tau)p_{\mu}, \\
\cC(\lambda; \tau; p)^{\bullet} & = & e^{\cC(\lambda; \tau; p)}.
\end{eqnarray*}

As pointed out in \cite{Mar-Vaf}, by comparing computations in
\cite{Mar-Vaf} with computations in \cite{Kat-Liu},
one obtains a conjectural formula for $\cC_{\mu}(\tau)$.
This formula is explicitly written down in \cite{Zho3}.
\begin{equation}\label{eqn:Mar-Vaf1}
\cC(\lambda; \tau; p)
=\sum_{n \geq 1} \frac{(-1)^{n-1}}{n}\sum_{\mu}
\left(\sum_{\cup_{i=1}^n \mu^i = \mu}
\prod_{i=1}^n \sum_{|\nu^i|=|\mu^i|} \frac{\chi_{\nu^i}(C(\mu^i))}{z_{\mu^i}}
e^{\sqrt{-1}(\tau+\frac{1}{2})\kappa_{\nu^i}\lambda/2} V_{\nu^i}(\lambda)
\right)p_\mu,
\end{equation}
\begin{equation}\label{eqn:Mar-Vaf2}
 \cC(\lambda;\tau; p)^{\bullet} = \sum_{|\mu| \geq 0}
\left(\sum_{|\nu|=|\mu|} \frac{\chi_{\nu}(C(\mu))}{z_{\mu}}
e^{\sqrt{-1}(\tau+\frac{1}{2})\kappa_{\nu}\lambda/2} V_{\nu}(\lambda)\right)
p_\mu,
\end{equation}
where
\begin{equation} \label{eqn:V}
\begin{split}
V_{\nu}(\lambda) = & \prod_{1 \leq a < b \leq l(\nu)}
\frac{\sin \left[(\nu_a - \nu_b + b - a)\lambda/2\right]}{\sin \left[(b-a)\lambda/2\right]} \\
& \cdot\frac{1}{\prod_{i=1}^{l(\nu)}\prod_{v=1}^{\nu_i} 2 \sin \left[(v-i+l(\nu))\lambda/2\right]}.
\end{split} \end{equation}

We now explain the notation on the right-hand sides of (\ref{eqn:Mar-Vaf1}) and
(\ref{eqn:Mar-Vaf2}). For a partition $\mu$ given by
$$\mu_1 \geq \mu_2 \geq \dots \geq \mu_{l(\mu)} > 0,$$
$\chi_{\mu}$ denotes the character of the irreducible representation of $S_d$
indexed by $\mu$, where $d=|\mu| = \sum_{i=1}^{l(\mu)} \mu_i$.
The number $\kappa_{\mu}$ is defined by
$$\kappa_{\mu} = |\mu| + \sum_i (\mu_i^2 - 2i\mu_i).$$
For each positive integer $i$,
$$m_i(\mu) = |\{j: \mu_j = i\}|.$$
Denote by $C(\nu)$ the conjugacy class of $S_d$ corresponding to the partition $\nu$,
and by $\chi_{\mu}(C(\nu))$ the value of the character $\chi_{\mu}$ on the conjugacy class $C(\nu)$.
Finally,
$$z_{\mu} =  \prod_j m_j(\mu)!j^{m_j(\mu)}. $$

In this paper, we will call (\ref{eqn:Mar-Vaf1}) the Mari\~{n}o-Vafa formula.

The third author proved in \cite{Zho1} some special cases of the
Mari\~{n}o-Vafa formula and found several interesting
applications. He showed in \cite{Zho4} that calculation of BPS
numbers in the local $\bP^2$ and $\bP^1\times\bP^1$ geometries can
be reduced to the Mari\~{n}o-Vafa formula, and proved in
\cite{Zho5} a special case of a conjecture by A. Iqbal \cite{Iqb} assuming
the Mari\~{n}o-Vafa formula.

We now describe our approach to the
Mari\~{n}o-Vafa formula (\ref{eqn:Mar-Vaf1}).
Denote the right-hand sides of (\ref{eqn:Mar-Vaf1})
and (\ref{eqn:Mar-Vaf2}) by $R(\lambda; \tau;p)$
and $R(\lambda;\tau;p)^{\bullet}$ respectively.
In \cite{Zho3}, the third author proved the following two equivalent cut-and-join
equations similar to the one satisfied by Hurwitz numbers \cite{Gou-Jac-Vai}, \cite{Li-Zha-Zhe},
\cite[Section 15.2]{Ion-Par2}.

\begin{Theorem}  \label{RCutJoin}
\begin{eqnarray}
&& \frac{\partial R}{\partial \tau}
= \frac{\sqrt{-1}\lambda}{2} \sum_{i, j\geq 1} \left(ijp_{i+j}\frac{\partial^2R}{\partial p_i\partial p_j}
+ ijp_{i+j}\frac{\partial R}{\partial p_i}\frac{\partial R}{\partial p_j}
+ (i+j)p_ip_j\frac{\partial R}{\partial p_{i+j}}\right), \label{eqn:CutJoin1} \\
&& \frac{\partial R^{\bullet}}{\partial \tau}
= \frac{\sqrt{-1}\lambda}{2} \sum_{i, j\geq 1} \left(ijp_{i+j}\frac{\partial^2R^{\bullet}}{\partial p_i\partial p_j}
+ (i+j)p_ip_j\frac{\partial R^{\bullet}}{\partial p_{i+j}}\right). \label{eqn:CutJoin2}
\end{eqnarray}
\end{Theorem}

Here is a crucial observation: One can rewrite (\ref{eqn:CutJoin2}) as a
sequence of systems of ordinary equations,
one for each positive integer $d$,
hence if $\cC(\lambda;\tau;p)^{\bullet}$ satisfies (\ref{eqn:CutJoin2}),
then it is determined by the initial value $\cC(\lambda; 0; p)^{\bullet}$.
To prove (\ref{eqn:Mar-Vaf1}) or (\ref{eqn:Mar-Vaf2}),
it suffices to prove the following two statements:
\begin{itemize}
\item[(a)] Equation (\ref{eqn:CutJoin1}) is satisfied by $\cC(\lambda; \tau; p)$.
\item[(b)] $\cC(\lambda; 0; p) = R(\lambda; 0; p)$.
\end{itemize}
Or equivalently,
\begin{itemize}
\item[(a)'] Equation (\ref{eqn:CutJoin2}) is satisfied by $\cC(\lambda; \tau; p)^\bullet$.
\item[(b)'] $\cC(\lambda; 0; p)^\bullet = R(\lambda; 0; p)^\bullet$.
\end{itemize}

It is shown in \cite{Zho3} that (b) holds. Therefore, the
Mari\~{n}o-Vafa formula (\ref{eqn:Mar-Vaf1})
follows from the following theorem.

\begin{Theorem}\label{CCutJoin}
\begin{eqnarray}
&& \frac{\partial \cC}{\partial \tau}
= \frac{\sqrt{-1}\lambda}{2} \sum_{i, j\geq 1} \left(ijp_{i+j}
\frac{\partial^2\cC}{\partial p_i\partial p_j}
+ ijp_{i+j}\frac{\partial \cC}{\partial p_i}\frac{\partial\cC}{\partial p_j}
+ (i+j)p_ip_j\frac{\partial \cC}{\partial p_{i+j}}\right) \label{eqn:CutJoin}
\end{eqnarray}
\end{Theorem}

The rest of the paper is organized
as follows. In Section \ref{initial},  we give a proof of the initial
condition (b). In Section \ref{combinatorics}, we give the proof of
Theorem \ref{RCutJoin} in \cite{Zho3}. In Section \ref{geometry}, we outline
the proof of Theorem \ref{CCutJoin} in \cite{LLZ}. The details we omit here
are straightforward calculations which will be given in \cite{LLZ}.
Complete lists of relevant references will be given in \cite{Zho3, LLZ}.

\section{Initial Condition}\label{initial}

The proof of the initial condition (b) needs the following two theorems.

\begin{theorem}
We have
\begin{eqnarray} \label{eqn:CInit}
&& \cC(\lambda; 0; p) = - \sum_{n > 0}  \frac{\sqrt{-1}^{n+1}p_n}{2n\sin (n\lambda/2)}.
\end{eqnarray}
\end{theorem}

\begin{proof}
When $l(\mu) > 1$,
we clearly have
$$\cC_{\mu}(\lambda;0) = 0.$$
When $\mu = (n)$ we have
\begin{eqnarray*}
\cC_{(n)} (\lambda; 0)
& = & - \sum_{g \geq 0}
\lambda^{2g-1} \sqrt{-1}^{n+1} \frac{\prod_{a=1}^{n-1} (n\cdot 0 + a)}{(n-1)!}
\int_{\Mbar_{g, 1}} \frac{\Lambda_g^{\vee}(1)\Lambda_g^{\vee}(0)\Lambda_g^{\vee}(-1)}
{1- n\psi_1} \\
& = & - \frac{\sqrt{-1}^{n+1}}{n} \sum_{g \geq 0}  (n\lambda)^{2g-1}
\int_{\Mbar_{g, 1}} \lambda_g \psi_1^{2g-2} \\
& = &  - \frac{\sqrt{-1}^{n+1}}{n^2\lambda} \cdot \frac{n\lambda/2}{\sin (n\lambda/2)} \\
& = & - \frac{\sqrt{-1}^{n+1}}{2n\sin (n\lambda/2)}.
\end{eqnarray*}
In the second equality we have used the Mumford's relations \cite[5.4]{Mum}:
$$\Lambda_g^\vee(1)\Lambda_g^{\vee}(-1)=(-1)^g.$$
In the third equality we have used \cite[Theorem 2]{Fab-Pan1}.
This proves (\ref{eqn:CInit}).
\end{proof}

\begin{theorem}
We have the following identity:
\begin{eqnarray}
&& \log\left(  \sum_{n \geq 0} \sum_{|\rho|=n}
\frac{e^{\frac{1}{4}\kappa_{\rho}\sqrt{-1}\lambda}}{\prod_{e \in \rho} 2\sin (h(e)\lambda/2)}
\frac{\chi_{\rho}(\eta)}{z_{\eta}} p_{\eta}\right)
= - \sum_{d \geq 1} \frac{\sqrt{-1}^{d+1}p_d}{2d\sin(d\lambda/2)}. \label{eqn:Key}
\end{eqnarray}
\end{theorem}

For a partition $\eta$,
$$n(\eta) = \sum_i (i-1)\eta_i
= \sum_i \begin{pmatrix} \eta'_i \\ 2\end{pmatrix}.$$
For any box $e \in \eta$,
denote by $h(e)$ its hook length.
Then
$$\sum_{x\in \eta} h(x) = n(\eta) + n(\eta') + |\eta|.$$

\begin{lemma}
Introducing formal variables $x_1, \dots, x_n, \dots$ such that
$$p_i(x_1, \dots, x_n, \dots) = x_1^i + \cdots + x_n^i + \cdots.$$
Then for for any positive integer $n$,
we have
\begin{eqnarray} \label{eqn:pg}
&& \sum_{n \geq 0} t^n\sum_{|\rho|=n}
\frac{q^{n(\rho)}}{\prod_{e \in \rho} (1 - q^{h(e)})}
\frac{\chi_{\rho}(\eta)}{z_{\rho}} p_{\eta}
= \frac{1}{\prod_{i, j} (1 -tx_iq^{j-1})}.
\end{eqnarray}
\end{lemma}

\begin{proof}
Recall the following facts about Schur polynomials:
\begin{eqnarray}
&& s_{\rho}(x) = \sum_{\eta} \frac{\chi_{\rho}(\eta)}{z_{\eta}} p_{\eta}(x), \label{eqn:s} \\
&& s_{\rho}(1, q, q^2, \dots)
= \frac{q^{n(\rho)}}{\prod_{e \in \rho} (1 - q^{h(e)})}, \\
&& \sum_{n \geq 0} t^n \sum_{|\rho|=n} s_{\rho}(x)s_{\rho}(y)
= \frac{1}{\prod_{i, j} (1 -tx_iy_j)}.
\end{eqnarray}
Combining these two identities,
one gets:
\begin{eqnarray*}
&& \sum_{n \geq 0} t^n \sum_{|\rho|=n}
\frac{q^{n(\rho)}}{\prod_{e \in \rho} (1 - q^{h(e)})}s_{\rho}(x)
= \frac{1}{\prod_{i,j} (1 - tx_iq^{j-1})}.
\end{eqnarray*}
The proof is completed by (\ref{eqn:s}).
\end{proof}

As a corollary we now prove (\ref{eqn:Key}).
First we need the following:

\begin{lemma}
For any partition $\rho$ we have
\begin{eqnarray}
&& \frac{1}{2} \sum_{e \in \rho} h(e) - n(\rho)
=  \frac{1}{4}\kappa_{\rho} + \frac{1}{2}|\rho|.
\end{eqnarray}
\end{lemma}

\begin{proof}
\begin{eqnarray*}
&& \frac{1}{2} \sum_{e \in \rho} h(e) - n(\rho) = \frac{1}{2} (n(\rho') - n(\rho) + |\rho|)  \\
& = & \frac{1}{2} (\sum_i \begin{pmatrix} \rho_i \\ 2 \end{pmatrix}
- \sum_i (i-1)\rho_i + |\rho|) \\
& = & \frac{1}{4} (\sum_i \rho_i(\rho_i - 1) - 2 \sum_i i\rho_i + 4|\rho|) \\
& = & \frac{1}{4}\kappa_{\rho} + \frac{1}{2}|\rho|.
\end{eqnarray*}
\end{proof}

Let $q = e^{-\sqrt{-1}\lambda}$, and $t = \sqrt{-1}q^{1/2}$,
then we have
\begin{eqnarray*}
&& \sum_{n\geq 0} t^n\sum_{|\rho|=n}
\frac{q^{n(\rho)}}{\prod_{e \in \rho} (1 - q^{h(e)})}
\frac{\chi_{\rho}(\eta)}{z_{\rho}} p_{\eta} \\
& = & \sum_{n \geq 0} \sqrt{-1}^nq^{n/2} \sum_{|\rho|=n}
\frac{q^{n(\rho)-\frac{1}{2}\sum_{e\in \rho} h(e)}}{\prod_{e \in \rho} (q^{-h(e)/2} - q^{h(e)/2})}
\frac{\chi_{\rho}(\eta)}{z_{\rho}} p_{\eta} \\
& = & \sum_{n \geq 0} \sqrt{-1}^nq^{n/2} \sum_{|\rho|=n}
\frac{q^{-\frac{1}{4}\kappa_{\rho} - \frac{1}{2}n}}{\prod_{e \in \rho} (q^{-h(e)/2} - q^{h(e)/2})}
\frac{\chi_{\rho}(\eta)}{z_{\rho}} p_{\eta} \\
& = & \sum_{n \geq 0} \sum_{|\rho|=n}
\frac{e^{\frac{1}{2}f_{\rho}(2)\sqrt{-1}\lambda}}{\prod_{e \in \rho} 2\sin (h(e)\lambda/2)}
\frac{\chi_{\rho}(\eta)}{z_{\rho}} p_{\eta}.
\end{eqnarray*}
Hence by (\ref{eqn:pg}),
\begin{eqnarray*}
&& \log\left(  \sum_{n \geq 0} \sum_{|\rho|=n}
\frac{e^{\frac{1}{4}\kappa_{\rho}\sqrt{-1}\lambda}}{\prod_{e \in \rho} 2\sin (h(e)\lambda/2)}
\frac{\chi_{\rho}(\eta)}{z_{\rho}} p_{\eta}\right) \\
& = & \log\left(  \sum_{n\geq 0} t^n\sum_{|\rho|=n}
\frac{q^{n(\rho)}}{\prod_{e \in \rho} (1 - q^{h(e)})}
\frac{\chi_{\rho}(\eta)}{z_{\rho}} p_{\eta}\right) \\
& = & \log  \frac{1}{\prod_{i, j} (1 -tx_iq^{j-1})}
= \sum_{i,j \geq 1} \sum_{d \geq 1} \frac{1}{d} t^d q^{d(j-1)} x_i^d \\
& = & \sum_{j \geq 1} \sum_{d \geq 1} \frac{1}{d} t^d q^{d(j-1)} p_d
= \sum_{d \geq 1} \frac{p_d}{d} \frac{t^d}{1 - q^d} \\
& = & - \sum_{d \geq 0} \frac{\sqrt{-1}^{d+1}p_d}{2d\sin(d\lambda/2)}.
\end{eqnarray*}

By (\ref{eqn:Key}) we have
\begin{eqnarray*}
R(\lambda;0;p) & = & \log \left(  \sum_{n \geq 0} \sum_{|\rho|=n}
\frac{\chi_{\rho}(\eta)}{z_{\eta}} e^{\frac{1}{4}\kappa_{\rho}\lambda} V_{\rho} p_{\eta}\right)
= - \sum_{d \geq 1} \frac{\sqrt{-1}^{d+1}p_d}{2d\sin(d\lambda/2)},
\end{eqnarray*}
where we have used the following identity proved in \cite{Zho3}:
$$V_{\mu} = \frac{1}{2^l \prod_{x \in \mu} \sin [h(x)\lambda/2]}.$$
Hence (b) is proved.

\section{Proof of Theorem \ref{RCutJoin} } \label{combinatorics}

Recall
$$c_{\mu} = \sum_{g \in C_{\mu}} g$$
lies in the center of the group algebra $\bC S_d$, hence it acts
as a scalar $f_{\nu}(\mu)$ on any irreducible representation $R_{\nu}$.
In other words,
let $\rho: S_d \to \End R_{\nu}$ be the representation indexed by $\nu$,
then
$$\sum_{g \in C(\mu)} \rho_{\nu}(g) = f_{\nu}(\mu)\id.$$
We need the following interpretation of $\kappa_{\nu}$ in terms
of character:
$$\kappa_{\nu} = 2 f_{\nu}(C(2)).$$
See e.g. \cite[(5)]{Oko}.
Here we use $C(2)$ to denote the class of transpositions.
We need the following result:

\begin{lemma}
Suppose $h \in S_d$ has cycle type $\mu$.
The product $C_{(2)} \cdot h$ is a sum of elements of $S_d$ whose type is either a cut or a join of $\mu$.
More precisely,
there are $ijm_i(\mu)m_j(\mu)$ elements obtained from $h$
by joining an $i$-cycle in $h$ to a $j$-cycle in $h$,
and there are $(i+j)m_{i+j}(\mu)$ elements obtained from $h$ by
cutting an $(i+j)$-cycle into an $i$-cycle and a $j$-cycle.
\end{lemma}

\begin{proof}
Denote by $[s_1, \dots, s_k]$ a $k$-cycle. Then
$$[s, t]\cdot [s, s_2, \dots, s_i, t, t_2, \dots, t_j] = [s, s_2, \dots, s_i][t, t_2, \dots, t_j],$$
i.e., an $(i+j)$-cycle is cut into an $i$-cycle and a $j$-cycle.
Conversely,
$$[s, t]\cdot [s, s_2, \dots, s_i][t, t_2, \dots, t_j] = [s, s_2, \dots, s_i, t, t_2, \dots, t_j],$$
i.e., an $i$-cycle and a $j$-cycle is joined to an $(i+j)$-cycle.
Hence for a permutation $h$ of type $\mu$,
$c_{(2)} \cdot h$ is a sum of all elements obtained from $h$ by either a cut or a join.
Fix a pair of $i$-cycle and $j$-cycle of $h$,
there are $i \cdot j$ different ways to join them to an $(i+j)$-cycle.
Taking into the account of $m_i(\mu)$ choices of $i$-cycles, and $m_j(\mu)$ choices of $j$-cycles,
we get
$$ijm_i(\mu)m_j(\mu)$$
different ways to obtain an element from $h$ by joining an $i$-cycle in $h$ to a $j$-cycle in $h$.
Similarly,
fix an $(i+j)$-cycle of $h$,
there are $i+j$ different ways to cut it into an $i$-cycle and a disjoint $j$-cycle in $h$.
And taking into account the number of $(i+j)$-cycles in $h$,
we get
$$(i+j)m_{i+j}(\mu)$$
different ways to obtain an element from $h$ by cutting an $(i+j)$-cycle
into an $i$-cycle and a $j$-cycle.
\end{proof}

Now we have for any $h \in S_d$ of cycle type $\mu$
\begin{eqnarray*}
&& \sum_{\mu}  f_{\nu}(2) \frac{\chi_{\nu}(\mu)}{z_{\mu}}p_{\mu} \\
& = & \sum_{\mu} \tr [f_{\nu}(2) \id \cdot \rho_{\nu}(h)]
\cdot  \prod_i \frac{p_i^{m_i(\mu)}}{i^{m_i(\mu)} m_i(\mu)!} \\
& = & \sum_{\mu} \tr [\sum_{g \in C(2)} \rho_{\nu}(g) \cdot \rho_{\nu}(h)]
\cdot  \prod_i \frac{p_i^{m_i(\mu)}}{i^{m_i(\mu)} m_i(\mu)!} \\
& = &  \sum_{\mu}\tr \rho_{\nu}(\sum_{g \in C(2)} g\cdot h)
\cdot  \prod_i \frac{p_i^{m_i(\mu)}}{i^{m_i(\mu)} m_i(\mu)!} \\
& = & \sum_{\mu} \left(\sum_{\eta \in J_{i,j}(\mu)} ijm_i(\mu)m_j(\mu) \chi_{\nu}(\eta)
+ \sum_{\eta \in C_{i,j}(\mu)} (i+j)m_{i+j}(\mu) \chi_{\nu}(\eta)\right)
\cdot  \prod_i \frac{p_i^{m_i(\mu)}}{i^{m_i(\mu)} m_i(\mu)!} \\
& = &  \left(ijp_{i+j}\frac{\partial}{\partial p_i}\frac{\partial}{\partial p_j}
+ (i+j)p_ip_j \frac{\partial}{\partial p_{i+j}}\right) \sum_{\eta} \frac{\chi_{\nu}(\eta)}{z_{\eta}} p_{\eta}.
\end{eqnarray*}
Here we have the following notations.
Let $\mu, \eta$ be two partitions, both represented by Young diagrams.
We write $\eta \in J_{i,j}(\mu)$ and $\mu \in C_{i,j}(\eta)$
if $\eta$ is obtained from $\mu$ by remove a row of length $i$ and a row
of length $j$,
then adding a row of length $i+j$.
It follows that
\begin{eqnarray*}
&& \frac{\partial R(\lambda;\tau;p)^{\bullet}}{\partial \tau} \\
& = & \frac{\sqrt{-1}\lambda}{2} \sum_{\mu, \nu}
\left( f_{\nu}(2) \frac{\chi_{\nu}(C(\mu))}{z_{\mu}} p_{\mu}\right)
e^{\sqrt{-1}(\tau+\frac{1}{2})\kappa_{\nu}\lambda/2}
V_{\nu}(\lambda)  \\
& = & \frac{\sqrt{-1}\lambda}{2}\left(ijp_{i+j}\frac{\partial}{\partial p_i}\frac{\partial}{\partial p_j}
+ (i+j)p_ip_j \frac{\partial}{\partial p_{i+j}}\right) \sum_{\eta} \frac{\chi_{\nu}(\eta)}{z_{\eta}} p_{\eta}
e^{\sqrt{-1}(\tau+\frac{1}{2})\kappa_{\nu}\lambda/2}
V_{\nu}(\lambda).
\end{eqnarray*}
This finishes the proof of Theorem 1.

\section{Proof of Theorem \ref{CCutJoin} } \label{geometry}

\subsection{Moduli space of relative morphisms}
We first describe the moduli space of stable relative morphisms to
$\bP^1$ used in \cite{Li-Son}. The moduli spaces of stable relative
morphisms are constructed by J. Li \cite{Li1}. The construction
in symplectic geometry was carried out independently by
Li-Ruan \cite{Li-Rua} and Ionel-Parker \cite{Ion-Par1, Ion-Par2}.

Let $\bP^1[m]$ be a chain of $m+1$ copies $\bP^1$,
such that the $i$-th copy is glued to the $(i+1)$-th copy at the point
$p_1^{(i)}$ for $i \leq m$. The first copy will be referred to as the
root component, and the other components will be called the bubble components
A point $p_1^{(m)}$ is fixed on the $(m+1)$-th component.
Denote by $\pi[m]: \bP^1[m] \to \bP^1$ the map which is identity on the root
component and contracts all the bubble components to $p_1^{(0)}$.

Let $\mu$ be a partition of $d>0$. Let
$\Mbar_{g, 0}(\bP^1, \mu)$ be the moduli
space of morphisms
$$f: (C,x_1,\dots,x_{l(\mu)}) \to \bP^1[m],$$
such that
\begin{enumerate}
\item $(C,x_1,\ldots,x_{l(\mu)})$ is a prestable curve
      of genus $g$ with $l(\mu)$ marked points.
\item $f^{-1}(p_1^{(m)})=\sum_{i=1}^{l(\mu)}$ as Cartier divisors,
      and $\deg(\pi[m]\circ f)=d$.
\item The preimage of each node in $\bP^1[m]$ consists of nodes of $C$.
      If $f(y)=p_1^{i}$ and $C_1$ and $C_2$ are two irreducible
      components of $C$ which intersects at $y$, then
      $f|_{C_1}$ and $f|_{C_2}$ has the same contact order to
      $p_1^{i}$ at $y$.
\item The automorphism group of $f$ is finite.
\end{enumerate}
Two such morphisms are isomorphic if they differ by an isomorphism
of the domain and an automorphism of the bubble components
of $\bP^1[m]$.
In particular, this defines the automorphism group in the stability
condition (4) above.

In \cite{Li1, Li2}, J. Li showed that $\Mbar_{g,0}(\bP^1,\mu)$ is a
separated, proper Deligne-Mumford stack with
a perfect obstruction theory of virtual
dimension
$$r=2g-2 + |\mu| + l(\mu),$$
so it has a virtual fundamental class of degree $r$.

\subsection{Torus action}

Consider the $\bC^*$-action
$$t \cdot [z^0:z^1] = [tz^0: z^1]$$
on $\bP^1$. It has two fixed points
$p_0= [0:1]$ and $p_1=[1:0]$.
This induces an action on $\bP^1[m]$ by the action on the root component
induced by the isomorphism to $\bP^1$, and the trivial actions on the
bubble components. This in turn induces an action on
$\Mbar_{g, 0}(\bP^1, \mu)$.

\subsection{The branch morphism}
There is a branch morphism
$$
\mathrm{Br}: \Mbar_{g,0}(\bP^1,\mu) \to
\mathrm{Sym}^r\bP^1\cong \bP^r.
$$
Note that $\bP^r$ can be identified with $\bP(H^0(\bP^1,\cO(r))$, and
the isomorphism
$$\bP(H^0(\bP^1,\cO(r))\cong\mathrm{Sym}^r\bP^1$$
is given by $[s]\mapsto \mathrm{div}(s)$.
The $\bC^*$ action on $\bP^1$ induces a
$\bC^*$ action on $H^0(\bP^1,\cO(r))$ by
$$
t\cdot (z^0)^k (z^1)^{r-k} = t^{-k}(z^0)^k (z^1)^{r-k}.
$$
So $\bC^*$ acts on $\bP^r$ by
$$
t\cdot [a_0,a_1,\ldots, a_r]= {[a_0,t^{-1}a_1,\ldots, t^{-r} a_r]},
$$
where $(a_0,a_1,\ldots, a_r)$ corresponds to
$\sum_{k=0}^r a_k (z_0)^k (z^1)^{r-k}\in H^0(\bP^1,\cO(r))$.
With this action, the branch morphism is $\bC^*$-equivariant.

\subsection{The Obstruction Bundle}

In \cite{Li-Son}, J. Li and Y. Song constructed
an obstruction bundle over the stratum where the target
is $\bP^1[0]=\bP^1$, and proposed an extension over the
entire $\Mbar_{g,0}(\bP^1,\mu)$.
Here we use a different extension which is equivalent
to the one used in \cite{Bry-Pan}.

Let $\pi[m]: \bP^1[m] \to \bP^1$ be the contraction to the
root component, and denote $\tilde{f} = \pi[m] \circ f$.
Dual to the obstruction space at a map $f: (C, x_1, \dots, x_{l(\mu)})
\to \bP^1[m]$,
consider the vector bundle $V$ with fiber at $f$ given by
$$H^1(C, \cO_C(-D)) \oplus H^1(C, \tilde{f}^*\cO_{\bP^1}(-1)),$$
where $D=x_1+\ldots+x_{l(\mu)}$.
It is a direct sum of two vector bundles $V_D$ and $V_{D_d}$.

Note that
$$H^0(C, \cO_C(-D))= H^0(C, \tilde{f}^*\cO_{\bP^1}(-1))=0,$$
so the ranks of $V_D$ and $V_{D_d}$ are, by Riemann-Roch,
$l(\mu)+g-1$ and $d+g-1$, respectively.

We lift $\bC^*$ action on $\Mbar_{g,0}(\bP^1,\mu)$ to
$V_D$ and $V_{D_d}$ as follows. The action on
$V_{D_d}$ comes from an action on $\cO_{\bP^1}(-1)\to\bP^1$
with weights $p$ and $p+1$ at the two fixed points $p_0$ and
$p_1$, respectively, where $p\in\bZ$. The fiber of
$V_D$ does not depend on the map $f$, so the fibers over
two points in the same orbit of the $\bC^*$ action can be
canonically identified. The action of $\lambda\in\bC^*$ on $V_D$
is multiplication by $\lambda^{-p-1}$.

\subsection{Functorial localization}
Let $T=\bC^*$. We will compute
$$
\mathrm{Br}_* e_T(V)=\sum_{l=0}^r a_l(p) H^l u^{r-l}.
$$
by virtual functorial localization \cite{LLY}.

Let $F(p,x)=\sum_{l=0}a_l(p)x^l$. We have
$$
\frac{f_k^*\mathrm{Br}_* e_T(V)}{e_T(T_{p_k}\bP^r)}=
\frac{F(p,k)}{(-1)^{r-k} k!(r-k)!}.
$$

By functorial localization, we have
$$
\int_{p_{r-k}}\frac{f_{r-k}^* \mathrm{Br}_* e_T(V)}{e_T(T_{p_{r-k}}\bP^r)}
=\sum_{F\subset\mathrm{Br}^{-1}(p_{r-k})}\int_{[F]^{vir}}\frac{e_T(V)}{
e_T(N_F^{vir}) }
$$
for $k=0,\ldots,r$, where $N_F$ is the virtual normal bundle of
the fixed loci $F$ in $\Mbar_{g,0}(\bP^1,\mu)$.
It is computed in \cite{LLZ} that
$$
\sum_{F\subset\mathrm{Br}^{-1}(p_{r-k}) }
\int_{[F]^{vir}}\frac{e_T(V)}{e_T(N_F^{vir}) }=
(p+1)^k J_{g,\mu}^k(p),
$$
where $J^k_{g,\mu}(p)$ is a degree $r-k$ polynomial in $p$, and
$$
J_{g,\mu}^k(-p-1)=(-1)^{d-l(\mu)+k} J_{g,\mu}^k(p).
$$
Moreover, we have
\begin{eqnarray*}
J_{g,\mu}^0(p)&=&\sqrt{-1}^{l(\mu)-|\mu|}\cC_{g,\mu}(p),\\
J_{g,\mu}^1(p)&=&\sqrt{-1}^{l(\mu)-|\mu|-1}\left(
\sum_{\nu\in J(\mu)}I_1(\nu) \cC_{g,\nu}(p)
+\sum_{\nu\in C(\mu)}I_2(\nu)\cC_{g,\nu}(p)\right.\\
& & \left. +\sum_{g_1+g_2=g}\sum_{\nu^1\cup \nu^2\in C(\mu)}
   I_3(\nu^1,\nu^2) \cC_{g_1,\nu^1}(p) \cC_{g_2,\nu^2}(p)\right).
\end{eqnarray*}
Here we use the notation in \cite{Li-Zha-Zhe}.
The set $J(\mu)$ (join) consists of partitions of $d$ of the form
$$
\nu=(\mu_1,\ldots,\hat{\mu}_i,\ldots,\hat{\mu}_j,\ldots,\mu_{l(\mu)},
     \mu_i+\mu_j)
$$
and the set $C(\mu)$ (cut) consists of partitions of $d$ of the form
$$
\nu=(\mu_1,\ldots,\hat{\mu}_i,\ldots,\mu_{l(\mu)},j,k)
$$
where $j+k=\mu_i$. The precise definitions of $I_1$, $I_2$, and $I_3$
can be found in \cite{Li-Zha-Zhe}. It follows from the definition
that (\ref{eqn:CutJoin}) in Theorem~\ref{CCutJoin} is equivalent to
$$
\frac{d}{d\tau}J^0_{g,\mu}(\tau)= - J^1_{g,\mu}(\tau).
$$

Since
\begin{eqnarray*}
&&F(p,x)\\
&=&\sum_{k=0}^r\frac{F(p,k)}{(-1)^{r-k} k!(r-k)!}
x(x-1)\cdots(x-k+1)(x-k-1)\cdots (x-r)\\
&=&\sum_{k=0}^r(p+1)^{r-k} J^{r-k}_{g,\mu}(p)
x(x-1)\cdots(x-k+1)(x-k-1)\cdots (x-r)\\
&=&\sum_{k=0}^r (p+1)^k J^k_{g,\mu}(p)
x(x-1)\cdots(x-(r-k-1))(x-(r-k+1))\cdots (x-r),
\end{eqnarray*}
therefore,
$$
\mathrm{Br}_* e_T(V)=
\sum_{k=0}^r(p+1)^k J^k_{g,\mu}(p)
H(H-u)\cdots(H-(r-k-1)u)(H-(r-k+1)u)\cdots (H-ru).
$$

\subsection{Final Calculations}

Let $\tau=-p-1$, then
\begin{eqnarray*}
\mathrm{Br}_* e_T(V)
&=&
\sum_{k=0}^r(-\tau)^k J^k_{g,\mu}(-\tau-1)
H(H-u)\cdots(H-(r-k-1)u)\\
&& \ \ \ \cdot(H-(r-k+1)u)\cdots (H-ru)\\
&=& \sum_{k=0}^r(-\tau)^k (-1)^{d-l(\mu)+k}J^k_{g,\mu}(\tau)
H(H-u)\cdots(H-(r-k-1)u)\\
&& \ \ \ \cdot(H-(r-k+1)u)\cdots (H-ru)
\end{eqnarray*}
Therefore,
$$
\mathrm{Br}_* e_T(V)
=(-1)^{d-l(\mu)}
 \sum_{k=0}^r \tau^k J^k_{g,\mu}(\tau)
H(H-u)\cdots(H-(r-k-1)u)(H-(r-k+1)u)\cdots (H-ru).
$$
For $i=0,\ldots,r-1$, we have
\begin{eqnarray*}
&  & H^iH(H-u)\cdots(H-(r-k-1)u)(H-(r-k+1)u)\cdots (H-ru)\\
& = &((H-(r-k)u)+(r-k)u)^i H(H-u)\cdots(H-(r-k-1)u)\\
 & & \cdot(H-(r-k+1)u)\cdots (H-ru)\\
& = &((r-k)u)^i H(H-u)\cdots(H-(r-k-1)u)(H-(r-k+1)u)\cdots (H-ru)
\end{eqnarray*}
since
$$
H(H-u)\cdots(H-ru)=0.
$$
Therefore,
$$
\int_{\bP^r}\mathrm{Br}_* e_T(V)H^i
=(-1)^{d-l(\mu)}u^i\sum_{k=0}^r(r-k)^i \tau^k J^k_{g,\mu}(\tau).
$$
Let $J^k_{g,\mu}(\tau)=\sum_{j=0}^{r-k}a_j^k \tau^j$. We have
$$
u^{-i}\int_{\bP^r}\mathrm{Br}_* e_T(V)H^i =
(-1)^{d-l(\mu)}\sum_{l=0}^r
\left(\sum_{j+k=l} (r-k)^i a_{j}^k\right)
\tau^l.
$$
Here is a crucial observation: as a polynomial in $\tau$,
$u^{-i}\int_{\bP^r}\mathrm{Br}_* e_T(V)H^i$
is of degree no more than $i$. Therefore,
$$
\sum_{j+k=l} (r-k)^i a_{j}^k=0
$$
for $0\leq i < l\leq r$. Now fix $l$ such that $1\leq l\leq r$.
We have
\begin{equation}\label{linear}
\sum_{k=0}^l (r-k)^i a_{l-k}^k =0,
\ \ \ \ 0\leq i<l,
\end{equation}
which is a system of $l$ linear equations of the
$l+1$ variables $\{a_{l-k}^k: k=0,\ldots,l\}$.

Both
$$
\{
(r-t)^i : i=0,\ldots,l-1
\}
$$
and
$$
\{
1, t, t(t-1), \ldots, t(t-1)\ldots(t-l+2)
\}
$$
are bases of the vector space
$$
\{ f(t)\in\bQ[t] :\deg(f)\leq l-1 \},
$$
so there exists an invertible $l\times l$ matrix
$(A_{ij})_{0\leq i,j\leq l-1}$ such that
$$
t(t-1)\cdots (t-i+1)=\sum_{j=0}^{l-1} A_{ij}(r-t)^j.
$$
In particular,
$$
k(k-1)\cdots (k-i+1)=\sum_{j=0}^l A_{ij}(r-k)^j.
$$
for $k=0,1,\ldots,l$, so (\ref{linear}) is equivalent to
$$
\sum_{k=0}^l k(k-1)\cdots (k-i+1) a_{l-k}^k =0,
\ \ \ \ 0\leq i<l,
$$
i.e.,
$$
\sum_{k=i}^l\frac{k!}{(k-i)!}a_{l-k}^k =0,
\ \ \ \ 0\leq i<l.
$$
The above equations can be rewritten as
$$
\left(
\begin{array}{ccccccc}
1 & 1      & \cdots & \cdots   & \cdots & \cdots & 1    \\
0 & 1!     & 2      & \cdots   & \cdots & \cdots & l   \\
0 & 0      & 2!     & 3\cdot 2 & \cdots & \cdots & l(l-1) \\
0 & 0      & 0      & 3!       & \cdots & \cdots & l(l-1)(l-2)\\
0 & \vdots & \vdots &\vdots    & \ddots & \vdots & \vdots\\
0 & \cdots & \cdots & \cdots   &0      & (l-1)! &l(l-1)\cdots 2
\end{array}\right)
\left(\begin{array}{c}a^0_l\\ a^1_{l-1}\\ \vdots\\ \vdots\\
   a^l_0 \end{array}\right)
= \left(\begin{array}{c}0\\\vdots\\ \vdots\\ \vdots\\ 0 \end{array}\right).
$$
The kernel is clearly one dimensional. One can check that the kernel
is given by
\begin{equation}\label{solution}
a_{l-k}^k=(-1)^k \frac{l!}{k!(l-k)!}a_l^0.
\end{equation}
Note that (\ref{solution}) for $l=1,\ldots,r$ is equivalent to
$$
J_{g,\mu}^k(\tau)=\frac{(-1)^k}{k!}\frac{d^k}{d\tau^k}J_{g,\mu}^0(\tau)
$$
for $k=0,\ldots,r$. In particular,
$$
J_{g,\mu}^1(\tau)=-\frac{d}{d\tau}J_{g,\mu}^0(\tau)
$$
which is equivalent to the cut-and-join equation (\ref{eqn:CutJoin})
in Theorem \ref{CCutJoin}.

\bigskip

{\bf Acknowledgements.}
We wish to thank Jun Li for explaining his work,
Jim Bryan, Bohui Chen, Tom Graber, Gang Liu for helpful conversations,
and Cumrun Vafa and Shing-Tung Yau for their interests in this work.
The first author wishes to thank the hospitality of IPAM where she
was a core participant of the Symplectic Geometry and Physics
Program and did most of her part of this work.
The second author is supported by an NSF grant.
The third author is partially supported by research
grants from NSFC and Tsinghua University.

\end{document}